\documentclass{amsart}

\usepackage{epsfig}
\usepackage{amscd}
\usepackage{amsmath, amssymb, amsthm}
\usepackage{graphics}
\usepackage{color}
\usepackage{dsfont}
\newtheorem*{main-theorem}{Main Theorem}
\newtheorem{proposition}{Proposition}[section]
\newtheorem{theorem}{Theorem}
\newtheorem{lemma}[proposition]{Lemma}
\newtheorem{corollary}[proposition]{Corollary}

\theoremstyle{definition}

\newtheorem{remark}[proposition]{Remark}

\numberwithin{equation}{section}

\def\tP{\widetilde{P}}
\def\tp{\tilde{p}}
\def\11{\mathbf{1}}

\def\reals{{\mathbb R}}
\def\cx{{\mathbb C}}

\def\Ci{{\mathcal C}^\infty}

\def\Im{\,\mathrm{Im}\,}

\def\WF{\mathrm{WF}_h\,}

\def\supp{\mathrm{supp}\,}

\def\id{\,\mathrm{id}\,}

\def\O{{\mathcal O}}

\def\neigh{\mathrm{neigh}\,}
\def\Op{\mathrm{Op}\,}

\def\phi{\varphi}
\def\half{{\frac{1}{2}}}
\def\dist{\text{dist}\,}

\def\be{\begin{eqnarray*}}
\def\ee{\end{eqnarray*}}
\def\ben{\begin{eqnarray}}
\def\een{\end{eqnarray}}

\def\lll{\left\langle}
\def\rrr{\right\rangle}

\def\L2R{L_{\text{Rest}}^2}

\def\HH{\mathcal{H}}

\def\tchi{\tilde{\chi}}
\def\L2c{L^2_{\text{comp}}}
\def\tphi{\tilde{\phi}}
\def\Dom{\text{Dom}\,}

\begin{document}
\title[Dispersive Estimates]{Dispersive Estimates for Manifolds 
with one Trapped Orbit}
\author{Hans Christianson}
\address{Department of Mathematics, MIT, 77 Mass. Ave., Cambridge, MA 02139-4307, USA}
\email{hans@math.mit.edu}

\begin{abstract}
For a large class of complete, non-compact Riemannian manifolds, $(M,g)$, with boundary, we prove high energy resolvent estimates in the case
where there is one trapped hyperbolic geodesic.  As an application, we have the following local smoothing estimate for
the Schr\"odinger propagator:
\be
\int_0^T \left\| \rho_s e^{it(\Delta_g-V)} u_0 \right\|_{H^{1/2 -
    \epsilon}(M)}^2 dt  \leq C_T \| u_0 \|_{L^2(M)}^2,
\ee
where $\rho_s(x) \in \Ci(M)$ satisfies $\rho_s = \lll \dist_g(x,x_0)
\rrr^{-s}$, $s> \half$, and $V \in \Ci(M)$, $0 \leq V \leq C$
satisfies $|\nabla V| \leq C \lll \dist(x,x_0) \rrr^{-1-\delta}$ for
some $\delta>0$.  From the
local smoothing estimate, we deduce a family of Strichartz-type
estimates, which are used to prove two well-posedness results for
the nonlinear Schr\"odinger equation.

As a second application, we prove the following sub-exponential local energy
decay estimate for solutions to the
wave equation when $\dim M = n \geq 3$ is odd and $M$ is equal to
$\reals^n$ outside a compact set:
\be
 \lefteqn{ \int_M \left|\psi \partial_t u \right|^2 + \left| \psi \nabla u
 \right|^2 dx } \\
& \leq & C e^{-t^{1/2} /C} \left( \|u(x,0)\|_{H^{1+\epsilon}(M)}^2 +
  \|D_tu(x,0)\|_{H^\epsilon(M)}^2 \right),
\ee
where $\psi \in \Ci(M)$, $\psi \equiv e^{-|x|^2}$ outside a compact set.

\end{abstract}
\maketitle


\section{Introduction}

In this note we show how the results of \cite{Ch, Ch2} on cutoff
resolvent estimates near closed hyperbolic orbits can be
combined with the non-trapping resolvent estimates in \cite{CPV} to
obtain resolvent bounds in the case of one trapped hyperbolic orbit
with a logarithmic loss.  As applications, we prove local smoothing
estimates for solutions to the linear Schr\"odinger equation (Theorem
\ref{Sch-thm}) and local
energy decay estimates for solutions to the linear wave equation
(Theorem \ref{Wave-thm}).  These theorems have direct applications to
the nonlinear Schr\"odinger and wave equations.

We prove the high-energy resolvent estimates for a much more general
class of manifolds, then specialize to the case of asymptotically
Euclidean manifolds for the applications.  The class of manifolds we consider
for the high-energy estimates are the same as those studied (in the non-trapping case) in \cite{CPV}.  More precisely, let $(M,g)$
be a connected Riemannian manifold, $M = X_0 \cup X$, where $X_0$ is a
compact, connected $n$-dimensional Riemannian manifold and
$X = [r_0, +\infty) \times S$, $r_0 \gg 1$, where $S$ is a compact, connected
  $(n-1)$-dimensional Riemannian manifold without boundary.  We assume
  $\partial X_0$ is compact and that $X$ and $X_0$ satisfy
\be
\partial X_0 = \partial M \cup \partial X,\,\,\, \partial M \cap
\partial X = \emptyset.
\ee
We assume the metric $g |_{X_0}$ is a $\Ci$ metric on
$\overline{X_0}$ and 
\be
g|_{X} = dr^2 + \sigma(r),
\ee
where $\sigma(r)$ is a family of smooth Riemannian metrics on $S$ depending
smoothly on $r$.  In local coordinates, the metric $\sigma(r)$ takes
the form
\be
\sigma(r) = \sum_{i,j = 1}^{n-1} g_{ij}(r, \theta) d \theta^i d
\theta^j,
\ee
and if we set $X_r = [r, +\infty) \times S$, we can identify $\partial
  X_r \simeq (S, \sigma(r))$.  Thus with $b = (\det g_{ij})^\half$ and
  $(g^{ij}) = (g_{ij})^{-1}$, we have
\be
\Delta_{\partial X_r} = -b^{-1}\sum_{i,j} \partial_{\theta_i}(b g^{ij}
\partial_{\theta_j}),
\ee
and
\be
\Delta_X = -b^{-1} \partial_r(b \partial_r) + \Delta_{\partial X_r}.
\ee
As in the introduction of \cite{CPV}, a calculation shows 
\be
b^\half \Delta_X b^{-\half} = \partial_r^2 + \Lambda_r + q(r, \theta),
\ee
with 
\be
\Lambda_r = - \sum_{i,j} \partial_{\theta_i}(g^{ij}
\partial_{\theta_j}),
\ee
and
\be
q(r,\theta) = (2b)^{-2} \left( \frac{ \partial b }{\partial r}
\right)^2 + (2b)^{-2} \sum_{i,j} \frac{\partial b }{\partial \theta_i}
\frac{\partial b}{ \partial \theta_j} g^{ij} + \half b
\Delta_X(b^{-1}).
\ee
We assume $q(r,\theta) = q_1(r,\theta) + q_2(r, \theta)$, where 
\be
&& |q_1(r, \theta)| \leq C, \,\,\, \left| \frac{\partial^k q_1}{\partial
  r^k} \right| \leq
Cr^{-k-\delta} \text{ for } k \geq 1, \text{ and } \\
&& \left| \frac{ \partial^{k'}q_2(r, \theta)}{\partial r^{k'}} \right|
\leq C r^{-k'-1-\delta} \text{ for } k' \geq 0,
\ee
for $C, \delta >0$.  Observe this is satisfied for Euclidean space
using a polar decomposition outside of a ball of radius $r_0$ (where
$b = r^{n-1} \alpha(\theta)$), and for asymptotically Euclidean or
conic 
manifolds.  Define $h \in \Ci([r_0, +\infty) \times T^*(\partial X_r))$ by 
\be
h(r, \theta, \xi) = \sum_{i,j} g^{ij}(r, \theta) \xi_i \xi_j,
\ee
and assume there is a constant $C>0$ such that for all $(\theta, \xi)
\in T^*(\partial X_r)$,
\be
-\frac{\partial h}{\partial r} (r, \theta, \xi) \geq \frac{C}{r} h(r,
\theta, \xi).
\ee

Let $-\Delta_g$ be the Laplace-Beltrami operator acting on functions,
with Dirichlet boundary conditions if $\partial M \neq \emptyset$, and
suppose $V
\in \Ci(M)$, $0 \leq V \leq C$ satisfies
\ben
\label{v-assump}
|\nabla V| \leq C \lll \dist(x,x_0) \rrr^{-1-\delta}
\een
for some $\delta>0$.

The operator $P: = -\Delta_g + V(x)$ is an unbounded operator
\be
P: \HH \to \HH,
\ee
where $\HH = L^2(M)$, with domain $H^2(M)$ or $H^2(M) \cap H_0^1(M)$
in the case $\partial M \neq \emptyset$.    
In order to study the operator
\be
P - \tau
\ee
for $\tau \in \cx$ in some neighbourhood of $\reals$, we use the following semiclassical rescaling for $-\Delta_g$.  For $z  \in   [E-\delta, E+\delta] +
i(-c_0h, c_0h)$ write
\be
\tau  = \frac{z}{h^2}.
\ee
Then 
\be
-\Delta_g +V(x) - \tau  & = & -\Delta_g +V(x) -
\frac{z}{h^2} \\
& = & \frac{1}{h^2} (-h^2 \Delta_g  +h^2 V(x) - z).
\ee

Now let $P(h) = -h^2 \Delta_g +h^2 V(x)$ be the self-adjoint semiclassical 
Schr\"odinger operator acting on $\HH$ with Dirichlet boundary
conditions if $\partial M \neq \emptyset$.  Let $p =
\sigma_h(P(h))$ be the semiclassical (Weyl) principal symbol of $P(h)$
(see \cite[Theorem D.3]{EvZw}).  We assume the Hamiltonian flow of
$H_p$ generates a single closed hyperbolic orbit 
$\gamma$ 
in the energy level $\{p = E\}$, $E>0$.  The assumption that $\gamma$
be hyperbolic means the linearization of the Poincar\'e map has no
eigenvalues on the unit circle (see \cite{Ch, Ch2} for definitions).
Let $\pi: T^*M \to M$ denote
the natural projection, and assume that the projected generalized
geodesic $\pi(\gamma)$ lies entirely within $U_0 \Subset U \Subset X_0$.  If
$\pi(\gamma)$ intersects $\partial M$, assume that the
intersection is transversal.  Assume further that
the geometry is non-trapping outside $U_0$.  That is, for every compact
subset $K \Subset M\setminus U_0$, there is a time $T(K)$ so that if $\eta(t)$ is
a generalized geodesic with $\eta(0) \in K$, there is a time $0 <
\tau \leq T(K)$ such that $\eta(\pm \tau) \in (M \setminus U_0) \setminus 
K$.

\subsection{The Main Results}

The following theorem is our local smoothing result for solutions to the
linear Schr\"odinger equation, and is a generalization of the results
in \cite{Bur2} and the references cited therein.  The Schr\"odinger propagator
$e^{it(\Delta_g -V(x))}$ is a unitary operator on $L^2(M)$, but this theorem says if we integrate in time, we gain some
regularity.
\begin{theorem}
\label{Sch-thm}
Suppose $(M,g)$ is a Riemannian manifold (with or without boundary)
which satisfies the above assumptions, $\gamma
\subset M$ is a closed hyperbolic geodesic, and $ -\Delta_g$ is the
Laplace-Beltrami operator (with Dirichlet boundary conditions if
$\partial M \neq \emptyset$).  Then for each $\epsilon>0$ and $T>0$, there is a constant $C$ such that
\ben
\label{Sch-thm-est}
\int_0^T \left\| \rho_s e^{it (\Delta_g - V(x)) } u_0 \right\|_{H^{1/2 -
    \epsilon}(M)}^2 dt  \leq C \| u_0 \|_{L^2(M)}^2,
\een
where $\rho_s \in \Ci(M)$ satisfies
\ben
\label{rho-def}
\rho_s(x) \equiv \lll d_g(x, x_0 )\rrr^{-s}
\een
for $x_0$ fixed and $x$ outside a compact set, and $V \in \Ci(M)$, $0 \leq V \leq C$
satisfies \eqref{v-assump}.
\end{theorem}

\begin{remark}
\label{rem-1a}
We will see that in some cases the weighted resolvent has no poles on the real
axis, and we can conclude the estimate
\eqref{Sch-thm-est} is global in time at the expense of replacing
$\rho_s$ with super-exponentially decreasing weights.  That is, in these cases
we have
\ben
\label{global-sm}
\int_0^\infty \left\| \psi e^{it (\Delta_g-V) } u_0 \right\|_{H^{1/2 -
    \epsilon}(M)}^2 dt  \leq C \| u_0 \|_{L^2(M)}^2,
\een
where $\psi \equiv \exp ( -\dist_g (x,x_0)^{2})$ outside a compact
set.  
This is the case, for example, if $g$ is an asymptotically Euclidean
scattering metric, $V \equiv 0$, and
$\partial M = \emptyset$ (see \cite[Theorem 3, \S 10]{Mel}).  It is
also the case if $(M,g)$ is equal to $\reals^n$ outside a compact set, $n
\geq 2$, and $V$ satisfies \eqref{v-assump-2} below (see \cite[Theorem 8, Ch.9]{Vai}).  See Remarks \ref{rem-1} and
\ref{rem-2}.
\end{remark}

As a second application, we study solutions to the linear wave equation on $(M,g)$:
\ben
\label{wave-eq-11}
\left\{ \begin{array}{l} (-D_t^2 - \Delta_g +V(x) ) u(x,t) = 0, \,\,\, (x,t) \in M
    \times [0,\infty) \\
u(x,0) = u_0 \in H^1(M), \,\, D_t u(x,0) = u_1 \in L^2(M), \end{array}
\right.
\een
where $-\Delta_g$ is the Dirichlet Laplace-Beltrami operator on
functions and $V\in \Ci(M)$ satisfies 
\ben
\label{v-assump-2}
\exp ( \dist_g (x,x_0)^{2}) V = o(1).
\een  
Let $\psi \in
\Ci(M)$ satisfy
\ben
\label{psi-def}
\psi \equiv \exp ( -\dist_g (x,x_0)^{2})
\een
for $x$ outside a compact set and $x_0$ fixed.  For $u$ satisfying
\eqref{wave-eq-11}, we define the {\it local energy}, $E_{\psi}(t)$, to be
\be
E_{\psi}(t) = \half \left( \left\| \psi \partial_t u
\right\|_{L^2(M)}^2 + \left\| \psi u \right\|_{H^{1}(M)}^2 \right)   .
\ee 

\begin{theorem}
\label{Wave-thm}
Suppose $(M,g)$ is equal to $\reals^n$ outside a compact set, $n = \dim M \geq
3$ is odd, and $\gamma \subset M$ is a hyperbolic trapped ray with no
other trapping.  Then for each
$\epsilon>0$ and each 
\be
u_0 &\in &\Ci_c(M) \cap H^{1+ \epsilon}(M), \,\,\, \text{and} \\
u_1 & \in & \Ci_c(M) \cap H^\epsilon (M),
\ee
there is a constant $C>0$ such that 
\ben
\label{loc-en-est}
E_{\psi}(t) \leq C e^{-t^{1/2}/C} \left( \|u_0\|_{H^{1+\epsilon}(M)}^2 +
  \|u_1\|_{H^\epsilon(M)}^2 \right).
\een
Here the constant $C$ depends only on $\epsilon>0$, $g$, $n$, $\psi$, and the
support of $u_0$ and $u_1$.
\end{theorem}

\begin{remark}
The estimate \eqref{loc-en-est} holds whenever the resolvent admits a
meromorphic extension to $\cx$ with no poles in a complex
neighbourhood of an interval
$[-C, C] \subset \reals$, which holds also,
for example, if $(M,g)$ is an exterior domain in $\reals^n$ with $n
\geq 3$ odd.  
\end{remark}

The problem of ``local smoothing'' estimates for the Schr\"odinger
equation has a long history.  The sharpest results to date are those
of Doi \cite{Doi} and Burq \cite{Bur2}.  Doi proved if $M$ is
asymptotically Euclidean, then one has the estimate
\ben
\label{doi-est}
\int_0^T \left\| \psi e^{it\Delta_g} u_0 \right\|_{H^{1/2}(M)}^2 dt  \leq C \| u_0 \|_{L^2(M)}^2
\een
for $\psi \in \Ci_c(M)$ if and only if there are no trapped sets.
Burq's paper showed if there is trapping due to the presence
of several convex obstacles in $\reals^n$ satisfying certain assumptions, then one has the estimate \eqref{doi-est}
with the $H^{1/2}$ norm replaced by $H^{1/2 - \epsilon}$ for $\epsilon
>0$.  

The estimates with the $\epsilon>0$ loss in trapping geometries corresponds to a
logarithmic loss in resolvent estimates for these geometries (see
Theorem \ref{res-thm-tau}).  With more care, one could replace the
$\epsilon>0$ loss in derivative with a logarithmic loss in derivative,
which may help in certain applications.  The proof of Theorem \ref{res-thm-tau} uses a semiclassical
reduction to consider an operator of the form
\be
P(h) -z = -h^2 \Delta_g  - z ,
\ee
with $z \in [E - \delta, E + \delta] + i (-c_0h, c_0h)$ for $E, \delta
>0$.  It is shown in \cite{Ch, Ch2} that for $A \in \Psi^{0,0}_h(M)$ with
sufficiently small wavefront set near $\gamma \subset \{ p = E \}$, if
$|\Im z | \leq c_0' h / \log(1/h)$, then 
\ben
\label{Ch-est}
\|  (P(h) - z) A u \|_{L^2} \geq C^{-1} \frac{h}{\log(1/h)} \| A
u \|_{L^2}.
\een
We will use the main results from \cite{CPV} and propagation of
singularities to extend this to an estimate on $M$.

As an application of Theorem \ref{Sch-thm} and the non-trapping
Strichartz estimates of
\cite{HTW}, we study the nonlinear Schr\"odinger equation
\ben
\label{nls}
\left\{ \begin{array}{cc}
i \partial_t u + (\Delta_g-V(x)) u = F(u) \,\, \text{on } I \times M; \\
u(0,x) = u_0(x),
\end{array} \right.
\een
where $I \subset \reals$ is an interval containing $0$ and $V \in
\Ci_c(M)$, $V \geq 0$.  Here the
nonlinearity $F$ satisfies
\be
F(u) = G'(|u|^2)u,
\ee
and $G: \reals \to \reals$ is at least $C^3$ and satisfies
\be
| G^{(k)}(t) | \leq C_k \langle t \rangle^{\beta - k },
\ee
for some $\beta \geq \half$.

In \S \ref{strichartz} we prove a family of Strichartz-type estimates  which will result in the following local
well-posedness proposition.  
(See \S \ref{strichartz} also for
comments on optimality.)  For the statement of the proposition, let $H_D^1(M)$ denote the domain of $(1- \Delta_g)^\half$ with
Dirichlet boundary conditions if $\partial M \neq \emptyset$, so that $H_D^1(M)
= H_0^1(M)$, and write $H_D^s(M)$ for the domain of $(1 -
\Delta_g)^{s/2}$ (with Dirichlet boundary conditions if $\partial M
\neq \emptyset$).

\begin{proposition}
\label{nls-lwp}
Suppose $(M,g$) satisfies the above assumptions, $V \in \Ci_c(M)$, and in
addition $M$ is asymptotically conic (as defined in \cite{HTW}).  Then for each 
\ben
\label{s-cond}
s > \frac{n}{2} - \frac{k}{\max \{2\beta -2, 2 \} }
\een
and each $u_0 \in H^s_D(M)$ there exists $p > \max \{ 2\beta -2, 2 \}$ and $0< T \leq
1$ such that \eqref{nls} has a unique solution
\ben
\label{u-soln}
u \in C ([-T, T]; H^s_D(M)) \cap L^p([-T, T]; L^\infty(M)).
\een
Here $k =1$ if $\partial M \neq \emptyset$ and $k = 2$ if $\partial M
= \emptyset$.  

Moreover, the map $u_0(x) \mapsto u(t,x) \in C([-T,T];H^s_D(M))$ is Lipschitz continuous on
bounded sets of $H^s_D(M)$, and if $\|u_0\|_{H^s_D}$ is bounded, $T$ is bounded
from below.
\end{proposition}

If we have $H^1$ energy conservation, Proposition \ref{nls-lwp}
implies $u$ extends to a global solution.

\begin{corollary}
\label{gwp-cor}
Suppose $(M,g)$ and $V$ satisfy the assumptions of Proposition
\ref{nls-lwp}, and assume $n \leq
2$.  If $G(r) \to + \infty$ as $r \to +
\infty$ then $u$ in \eqref{u-soln} extends to a solution
\be
u \in C ((-\infty, \infty); H^1_D(M)) \cap L^p((-\infty, \infty);
L^\infty(M)). 
\ee

If $\partial M = \emptyset$, $n \leq 3$, $\beta <3$, and $G(r) \to +
\infty$ as $r \to + \infty$, then the same conclusion holds.
\end{corollary}
\begin{remark}
In particular, the cubic defocusing non-linear Schr\"odinger equation
is globally well-posed.  Observe also that three spatial dimensions is
the smallest dimension in which the periodic orbit $\gamma$ can have a
Poincar\'e map whose linearization possesses complex eigenvalues.
\end{remark}

Local energy decay for solutions to the linear wave equation has also enjoyed a long
history.  Studied in non-trapping exterior domains by Morawetz
\cite{Mor}, Morawetz-Phillips \cite{MoPh}, and
Morawetz-Ralston-Strauss \cite{MRS}, and generalized by, for example
Vodev \cite{Vod}, it is well-known (see \cite{Ral})
that when there {\it are} trapped rays, one cannot expect exponential decay of the energy with
no loss in regularity.  Metcalfe-Sogge \cite{MeSo}
have recently shown that if there are trapped hyperbolic rays {\it
  and} sub-exponential energy decay with loss in derivative, then one has
long-time existence for certain classes of quasi-linear wave equations
in $\reals^n$.  Theorem \ref{Wave-thm} says this always happens with
one trapped hyperbolic orbit.  Specifically, suppose $M = \reals^n \setminus U$ for
$U \Subset \reals^n$, $-\Delta$ is the Dirichlet Laplacian, 
\be
Q (z, w) \in \Ci(\cx^n \times \cx^{n^2})
\ee
satisfies
\be
&& \text{ i) } Q \text{ is linear in } w, \\
&& \text{ii) For each }w,\,\, Q( \cdot, w) \text{ is a symmetric quadratic form,}
\ee
and consider the following initial value problem:
\ben
\label{nlw}
\left\{ \begin{array}{c}
(-D_t^2 - \Delta) u = Q(Du, D^2u) \text{ on } M \times [0,\infty) ,\\
u(x,0) = u_0, \,\, D_t u(x,0) = u_1.
\end{array} \right.
\een
The following Proposition then follows directly from \cite[Theorem 1.1]{MeSo} in
dimension $n=3$ and \cite[Theorem 1.1]{MeSo1} in dimensions $n\geq 5$.
\begin{proposition}
\label{meso-prop}
Suppose $(u_0 , u_1) \in (\Ci (\reals^n \setminus U))^2$, $n\geq 3$ odd, satisfy the
compatibility condition from \cite[\S 1]{MeSo}, and $\gamma \subset (\reals^n
\setminus U)$ is a trapped hyperbolic geodesic, with no other
trapping.  Assume further that if $n=3$, the null condition
\cite[(1.9), (1.10)]{MeSo} holds.  Then there exist $\epsilon_0>0$ and $N>0$ such that for every
$\epsilon \leq \epsilon_0$, if
\be
\sum_{|\alpha| \leq N} \| \lll x \rrr^{|\alpha|} \partial_x^\alpha u_0
\|_{L^2} + \sum_{|\alpha| \leq N-1} \| \lll x
\rrr^{|\alpha|+1} \partial_x^\alpha u_1 \|_{L^2} \leq \epsilon,
\ee
then \eqref{nlw} has a unique solution $u \in \Ci ([0,\infty) \times
\reals^n \setminus U)$.
\end{proposition}

{\bf Acknowledgments.}  
The author would like to thank Maciej Zworski for suggesting
  the local smoothing problem, as well as providing much help and support during the
  writing of this paper.  He would also like to thank Jason Metcalfe
  for helpful comments and suggesting Proposition \ref{meso-prop},
  Daniel Tataru for helpful discussions about optimality of Strichartz
  estimates, and Nicolas Burq for suggesting Proposition
  \ref{str-prop-5}.  The bulk of this paper was written while the
  author was a graduate student at UC-Berkeley, so he would like to
  thank the Mathematics Department at UC-Berkeley for their support.  Finally, he would like to thank the anonymous
  referee whose many comments and suggestions helped improve the exposition.


\section{Resolvent estimates}

Let 
\be
 ( P  - \tau )^{-1} = ( - \Delta_g  +V(x) - \tau )^{-1}
\ee
be the classical resolvent.  In this note we use the notation $\tau$ for the
unsquared spectral parameter and $\lambda^2 = \tau$ for the squared parameter.  It will be convenient to use the
lower half-plane as the physical half-plane.  The proof of Theorem \ref{Sch-thm} relies on the
weighted resolvent estimates of the following Theorem.

\begin{theorem}
\label{res-thm-tau}
Suppose $(M,g)$ satisfies all of the assumptions above.  Then for each
$\epsilon>0$ sufficiently small and each $s > \half$ there is
a constant $C>0$ such that
\ben
\label{res-tau-est}
\left\| \rho_s (P - (\tau \pm i
\epsilon ) )^{-1} \rho_s
\right\|_{\HH \to \HH } \leq
C \frac{ \log ( 2 + |\tau|) }{\lll \tau \rrr^{1/2}}, \,\,\, \tau \in \reals.
\een
\end{theorem}

\begin{remark}
\label{rem-1}
To prove \eqref{res-tau-est} is uniform in $\epsilon>0$,
it 
suffices by Proposition \ref{resolvent-prop} to prove the uniformity for $|\tau
| \leq C$ for some $C>0$.  This is the case if there are no embedded
eigenvalues in $\reals$.  This happens, for example, if $g$ is an
asymptotically Euclidean scattering metric and $\partial M =
\emptyset$, or if $(M,g)$ is equal to $\reals^n$ outside a compact
set.  In the latter case, for $\psi$ satisfying \eqref{psi-def}, 
\ben
\label{res-lam}
\psi (P - \lambda^2 )^{-1} \psi
\een
continues meromorphically to 
\be
\lambda  \in \left\{ \begin{array}{l} \cx, \,\,\, n \text{ odd}, \\ (\cx
    \setminus \{0 \} )^*, \,\,\, n \text{ even}, \end{array} \right.
\ee
where $( \cx \setminus \{ 0 \} )^*$ is the logarithmic Riemann
surface.  If, in addition, $V (x)$ satisfies \eqref{v-assump-2}, there
is no pole at $\lambda=0$, and \eqref{res-tau-est} is
uniform in $\epsilon>0$ (see \cite[Theorem 8, Ch. 9]{Vai}).

The contours we will be using are pictured in Figures \ref{spec-fig}
and \ref{spec2-fig}.  For details on the meromorphic continuation, see, for example, \cite{Sjo}.
\end{remark}

\begin{figure}
\hfill
\begin{minipage}[t]{.45\textwidth}
\centerline{\input{spec}}
\caption{\label{spec-fig} The curve $\tau - i \epsilon$ in the $z \in
  \cx$ plane.}
\end{minipage}
\hfill
\begin{minipage}[t]{.45\textwidth}
\centerline{
\input{spec2}}
\caption{\label{spec2-fig} The same curve in the $z^\half$ plane.}
\end{minipage}
\hfill
\end{figure}

To prove Theorem \ref{res-thm-tau} in general, we observe 
\be
\| \rho_s   (P - (\tau \pm i
\epsilon ) )^{-1} \rho_s \|_{\HH \to \HH} \leq \frac{C}{\epsilon}.
\ee
Using this estimate for $|\tau| \leq C$, we need only show \eqref{res-tau-est}
for $|\tau|$ large, which is Corollary \ref{res-cor}.

It is well known (see, for
example, \cite{BP}) that for $R>0$ sufficiently large, one can
construct a metric $\tilde{g}$ with no trapped geodesics so that
$\tilde{g}|_{X_R} = g|_{X_R}$.  Let
$\chi_s \in \Ci(M)$, $\supp \chi_s \subset X_{R+1}$, and $\chi_s(x) \equiv d_g(x,x_0)^{-s}$ for fixed $x_0$
and $x$ outside a compact set.  
If $\Delta_0$ is the Laplace-Beltrami
operator associated to $\tilde{g}$, we have
\be
 \Delta_g \chi_s =  \Delta_0 \chi_s,
\ee
 
Proposition $2.3$ and the Remark immediately following from \cite{CPV} show if $s
>1/2$, and $\chi \in \Ci(M)$, $\chi \equiv 1$ on $\supp \chi_s$ and
$\supp \chi \subset X_R$, then 
\ben
\label{CPV-est}
\left\| \chi_{-s} (P(h) - E \pm i \epsilon) \chi u \right\|_{L^2(M)} \geq C h \|\chi_s u \|_{H^1_h (M )}
\een
for $h>0$ sufficiently small.  Here $H^1_h(M)$ is the semiclassical
Sobolev space equipped with the norm
\be
\|u \|_{H^1_h(V)}^2 = \|u \|_{L^2(V)}^2 + \|h \nabla u \|_{L^2(V)}^2.
\ee
We prove the presence of $\gamma$
forces a weaker estimate.  
\begin{proposition}
\label{resolvent-prop}
Let $(M,g)$ satisfy the above assumptions.  
For each $\rho_s \in \Ci(M)$ satisfying \eqref{rho-def} there exist constants
$C, h_0 >0$ such
that for $0 < h \leq h_0$ 
\ben
\label{res-est}
\left\| \rho_s (P(h) - (E \pm i \epsilon))^{-1} \rho_s \right\|_{\HH \to \HH} \leq C
h^{-1} \log(1/h),
\een
uniformly in $\epsilon>0$.
\end{proposition}

We remark that an estimate similar to \eqref{Ch-est} was obtained in
\cite{BZ} under some more assumptions, and in that work the authors
implicitly suggested a result such as Proposition \ref{resolvent-prop}
should be possible.  

From Proposition \ref{resolvent-prop} we will be able to deduce the
following Corollary by rescaling.  We state a version both for $\tau$
and for $\lambda$.

\begin{corollary}
\label{res-cor}
Let $(M,g)$ satisfy the above assumptions.  For each
$\rho_s \in \Ci(M)$ satisfying \eqref{rho-def}, there exists a constant $C$ such that
\be
\left\| \rho_s (-\Delta_g + V(x) - \tau)^{-1}
  \rho_s \right\|_{\HH \to \HH} \leq C \frac{\log (2 + |\tau
  |)}{\lll \tau \rrr^{1/2}},
\ee
for $|\tau| \geq C$ and 
\be
| \Im \tau | \leq  \frac{\lll \tau \rrr^{1/2}}{C \log (2 + |\tau
  |)}.
\ee
Furthermore,
\be
\left\| \rho_s (-\Delta_g + V(x) - \lambda^2)^{-1}
  \rho_s \right\|_{\HH \to \HH} \leq C' \frac{\log (2 + |\lambda
  |)}{\lll \lambda \rrr},
\ee
for $| \lambda | \geq C'$ and 
\be
| \Im \lambda | \leq \frac{1}{C' \log (2 + | \lambda|)}.
\ee
\end{corollary}

\begin{proof}[Proof of Proposition \ref{resolvent-prop}]
Observe for $\pm \Im z \geq c_0h/ \log (1/h)$, \eqref{res-est} holds automatically so we need only prove the
Proposition for $| \Im z| \leq c_0 h/ \log (1/h)$ for some small constant
$c_0>0$.  Let $z \in [E-\delta, E + \delta] - i( c_0 h/ \log (1/h),0)$, $\delta>0$, for the remainder of the proof.

The idea of the proof will be to glue two cutoff resolvent estimates together
and control the interaction terms by propagation of singularities, and
then replace the cutoffs with $\rho_s$, again controlling the errors
with propagation of singularities.  There are $4$ main steps.

\vspace{.25cm}

\noindent {\bf Step 1: Select cutoffs.}

Recall we have defined $X_r = [r, + \infty) \times S$
in the introduction, chosen $R_0>0$
sufficiently large so that we can construct $\tP(h) = \Op(\tp)$ which agrees
with $P(h)$ on $X_{R_0}$, and the Hamiltonian flow of $\tp$ is
globally non-trapping.  Choose
$\psi \in \Ci_c(M)$, $0 \leq \psi \leq 1$, $\psi
\equiv 1$ on $M \setminus X_{R_0}$, $\psi \equiv 0$ on $X_{R_0 +1}$, 
and select $\chi_0, \chi_1 \in \Ci_c(M)$, $0 \leq \chi_0 \leq 1$, $\chi_0 \equiv 1$
near $\gamma$ with small support, $\chi_1 = \psi - \chi$.  


In order to control the interaction (commutator) terms, we will add a complex absorption potential to $P(h)-z$ which
is supported away from the above cutoffs, which will control the
interactions through propagation of singularities.  Choose $R_j$, $j=1,\ldots,7$,
\be
R_0 +1=: R_1, \,\,\, R_j < R_{j+1} < \infty,
\ee
and let
\be
A_{R_{j_1},R_{j_2}} = X_{R_{j_1}} \setminus X_{R_{j_2}}
\ee
be the annulus with inner radius $R_{j_1}$ and outer radius
$R_{j_2}$.  We will fix the distances between the $R_j$s at the end of
the proof.

Choose $a \in \Ci_c(M)$, $a \psi \equiv 0$, 
\be
a \equiv 1 \text{ on } A_{R_2, R_5}, \,\,\, \supp a \subset A_{R_1, R_6}
\ee
and choose $\psi_1, \psi_2 \in \Ci_c(M)$ satisfying $\psi_1  =
\psi_2^2$ and
\be
 \supp \psi_2 \subset A_{R_2,R_5}, \,\, \psi_2 \equiv 1 \text{ on }
A_{R_3, R_4}.
\ee
Set $Q(z) = P(h) - z - i C_1ha$ for a constant $C_1>0$ to be chosen later in
the proof.  

Recall 
\be
\gamma \Subset U_0 \Subset U \Subset X_0,
\ee
and choose $\tchi \in \Ci_c(M)$ satisfying 
\be
\tchi & \equiv & 1 \text{ on } M \setminus X_{R_6} \setminus U \text{
  and } \\
\supp \tchi & \subset & M \setminus X_{R_7} \setminus U_0.
\ee
Without loss of generality, we assume $\rho_s$ from the statement of
the Proposition satisfies \eqref{rho-def} and $\rho_s \equiv 1$ on $M
\setminus X_{R_7}.$  These cutoffs are shown pictorially in Figure \ref{fig-cuts}.


\begin{figure}
\input{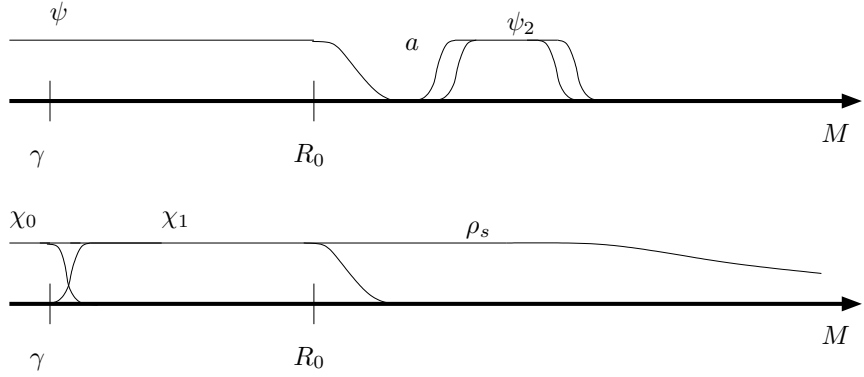}
\caption{\label{fig-cuts} The manifold $M$ with various cutoff
  functions employed in the proof of Proposition
  \ref{resolvent-prop}.}
\end{figure}

We will also employ an energy cutoff, to separate the characteristic
variety of $p-E$ from the elliptic sets.  Choose $\phi_1 \in \Ci_c(
\reals)$, 
\be
\phi_1(t) \equiv 1 \text{ on } \{ |t| \leq \alpha /2 \}, \,\,\,
\phi_1(t) \equiv 0 \text{ on } \{ |t| \geq \alpha \}.
\ee
Set
\be
\phi(x,\xi) = \phi_1 ( p(x, \xi) - E ),
\ee
and observe since $\phi$ is a function of the principal symbol of $P(h) - z$ 
and we are using the Weyl calculus,
\ben
\label{P-phi-comm}
[P(h) -z, \phi^w] = \O(h^3).
\een

\vspace{.25cm}

\noindent {\bf Step 2: Microlocalization.}

We will bound $\| \psi u \|$ from above, where unless explicitly
noted, $\| \cdot \| = \| \cdot \|_{\HH}$.  To do this, calculate
\be
\| \psi u \|   \leq \| \chi_0  u \| + \| \chi_1  u \| =: A + B.
\ee

For $B$ we cutoff in energy to apply \cite[Theorem 1,
Corollary 8]{Ch} and the generalizations from \cite{Ch2} for $c_0>0$
sufficiently small:
\ben
B & \leq & \| \phi^w \chi_0  u \| + \| (1- \phi)^w \chi_0  u \| \nonumber \\ 
& \leq & \frac{C \log (1/h)}{h} \left\|  (P(h) -z) \phi^w \chi_0  u
\right\| + \frac{C}{ \alpha} \left\| (P(h) - z) (1-\phi)^w \chi_0  u \right\|
\nonumber \\
& \leq & \frac{C \log (1/h)}{h} \left\| \phi^w (P(h) -z ) \chi_0  u
\right\| + C h^2 \log (1/h) \| \chi_0  u \| \label{st2-eq2-1} \\
&& + \frac{C}{ \alpha} \left\| (1 - \phi)^w (P(h) -z ) \chi_0  u \right\|
+ \frac{C}{ \alpha} h^3 \left\| \chi_0  u \right\| \label{st2-eq2-1a} \\
&\leq & \frac{C \log (1/h)}{h} \left( \left\| Q(z)  u
\right\| + \left\| [P(h), \chi_0 ]  u \right\| \right)  + C h^2 \log
(1/h) \| \chi_0  u \|, \label{st2-eq2-2} 
\een
since $\chi a = 0$.  Here in (\ref{st2-eq2-1}-\ref{st2-eq2-1a}) we have used
\eqref{P-phi-comm}. 
\vspace{.25cm}

To estimate $A$ and the commutator term in \eqref{st2-eq2-2} we will
need the lemmas in Step 3.

\noindent {\bf Step 3: Two Lemmas.}

The first Lemma is a refinement of the standard propagation of singularities result.

\begin{lemma}
\label{wf-lemma-2}
Let $\widetilde{V}_1, \widetilde{V}_2
\Subset M$, and for $j = 1,2$ let $V_j \Subset T^*M$,
\be
V_j := \{ (x, \xi) \in T^*M : x \in \widetilde{V}_j, \,\, |p(x,\xi)
-E| \leq \alpha \}, 
\ee
for some $\alpha>0$.  
Suppose the $\widetilde{V}_j$ satisfy $\dist_g(\widetilde{V}_1, \widetilde{V}_2) =
L,$ and assume 
\begin{eqnarray}
\label{dyn-assumption-10}
\left\{ \begin{array}{l}
\exists C_1,C_2 >0 \text{ such that }
\forall \rho \text{ in a neighbourhood of } V_1, \\
\exp( t H_p)(\rho) \in  V_2 \text{ for } \\
\sqrt{E}(L + C_1) \leq  t
\leq \sqrt{E}(L + C_1 + C_2).
\end{array} \right.
\end{eqnarray}
Suppose $A\in \Psi_h^{0,0}$ is microlocally equal to $1$ in $V_2  $.
If $B \in \Psi_h^{0,0}$ and $\WF (B) \subset V_1$, then there exists a
constant $C>0$ depending only on $C_1, C_2$ such that
\begin{eqnarray*}
\left\| Bu \right\| 
& \leq & C L h^{-1} \|B \|_{\HH \to \HH}  \left\| (P(h)-z)u
\right\| + 2 (E + \alpha)^{3/4} \frac{(C_1+1)}{\sqrt{C_2}} \|B\|_{\HH
  \to \HH} \| Au \|  \\
&& \quad + \O (h) \|\widetilde{B} u\|,
\end{eqnarray*}
where 
\be
\widetilde{B} \equiv 1 \text{ on } \cup_{0 \leq t \leq \sqrt{E}(L + C_1 + C_2)} \exp( t H_p)(
\WF B).
\ee
\end{lemma}

\begin{remark}
Observe Lemma \ref{wf-lemma-2} is a statement about the principal
symbol $p-z$, and hence applies also to $Q(z)$, and the difference is
$\O(h^\infty) \| \widetilde{B} u \|$.
\end{remark}

\begin{proof}
Let $G = \Op^w(g) \in
\Psi_h^{0,0}$ be a self-adjoint operator to be determined later in the proof and
calculate
\ben
\frac{L^2}{2h^2} \| G (P-z) u \|^2 + \frac{h^2}{2L^2}  \| G u \|^2 & \geq &
\Im \lll G (P-z) u, G u \rrr \nonumber \\ 
&  = & \Im \lll [G,P]u,Gu \rrr \nonumber \\
& \geq & \frac{h}{2} \lll \Op^w( \{ p, g^2 \}) u,u \rrr - \O(h^3) \|
\widetilde{G} u \|^2, \label{eqn-99}
\een
where $\widetilde{G} \equiv 1$ on $\WF G$.  
Hence 
\be
\frac{L^2}{2h^2} \| G (P-z) u \|^2 \geq \frac{h}{2} \lll \Op^w( \alpha) u,u \rrr
- \O(h^3) \|
\widetilde{G} u \|^2,
\ee
where
\be
\alpha(x,\xi) = \{ p, g^2\} - \frac{h}{T^2} g^2.
\ee
Choose $\phi_0 \in \Ci (M)$, $0 \leq \phi_0 \leq 1$, satisfying 
\be
\phi_0 & = & \phi_1^2, \text{ for } \phi_1 \in \Ci(M), \\
\phi_0 & \equiv &  1 \text{ on } \widetilde{V}_1,\\
\supp |\nabla \phi_0 | & \subset & \widetilde{V}_2, \text{ and} \\
| \nabla \phi_0 | & \leq & \frac{2}{C_2}.
\ee
Choose also $\phi \in \Ci( T^*M)$, $\phi = \phi(p(x,\xi) - E)$ so that
$\phi^2  \equiv 1$ on $V_1 \cup V_2$.
According to \eqref{dyn-assumption-10}, we can find a
non-characteristic hypersurface $\Sigma$ near $V_1$ so that 
\be
V_1 \cup V_2 \Subset \bigcup_{0 \leq t \leq \sqrt{E}(L + C_1 + C_2)}
\exp(t H_p)(\Sigma) =: \widetilde{\Sigma}.
\ee
Choose $f \in \Ci_c( \Sigma)$, $0 \leq f \leq 1$ so that $V_1$ and
$V_2$ are contained also in the flowout of $\{f = 1\}$, and choose
$\chi_0 \in \Ci_c(T^*M)$, $0 \leq \chi_0 \leq 1$, satisfying $\chi_0
\equiv 1$ on $V_1$ and $\chi_0 \equiv 0$ outside a neighbourhood of
$V_1$.  Let $q,a_0 \in \Ci(T^*M)$ be the solutions to
\be
&& H_p q = \chi_0, \,\,\, q|_\Sigma = f, \\
&& H_p a_0 = 1, \,\,\, a|_\Sigma = 0.
\ee
Observe $q$ satisfies
\be
&& 1 \leq q \leq \sqrt{E + \alpha} C_1 \text{ on } V_1; \\
&& |q| \leq \sqrt{E + \alpha} (C_1 +1 ) \text{ on }
\widetilde{\Sigma},
\ee
if $\supp \chi_0$ is sufficiently small.  In addition, $a_0$ satisfies
\be
0 \leq a_0 \leq \sqrt{E + \alpha} (L + C_1 + C_2) \text{ on }
\widetilde{\Sigma}.
\ee
Set
\be
g^2 = \phi^2( p - E) \phi_0(x) q^2 \exp (2h a_0/L^2), 
\ee
so that with this choice of $g^2$,
\ben
\alpha & = & \{ p, g^2\} - \frac{h}{T^2} g^2 \nonumber \\
& = & 2q \{ p, q \} e^{2h a_0/L^2} + 2\frac{h}{T^2} g^2 + q^2 e^{2h
  a_0/L^2} \phi^2\{p, \phi_0 \} - \frac{h}{T^2} g^2 \nonumber \\
& \geq &  2q \{ p, q \} e^{2h a_0/L^2} + \frac{h}{T^2} g^2 - 2
(E+\alpha)^{3/2} \frac{(C_1+1)^2}{C_2}  \label{eqn-100} 
\een
Combining \eqref{eqn-99} with \eqref{eqn-100} gives the lemma.
\end{proof}

The second Lemma will follow from Lemma \ref{wf-lemma-2} and indicates
how to control the interaction terms of the form $\| [P, \chi ] u \|$.

\begin{lemma}
Suppose $\chi \in \Ci_c(M)$ satisfies 
\be
\supp \chi \subset M \setminus X_{R_0}, \text{ and } \nabla \chi
\equiv 0 \text{ near } \gamma.
\ee
Then
\be
h^{-1} \left\| [P(h), \chi ] u \right\| \leq C R_7 h^{-1} \| Q(z)  u \| + \O(
h ) \| \tchi u \|,
\ee
where $\tchi \in \Ci_c(M)$ was selected in Step 1.
\end{lemma}

\begin{proof}
We first microlocalize using $\phi^w$ as in Step 2.  Observe
$[P(h),\chi] = h A(x,hD) $, where $A(x, hD)$ is a first order
semiclassical differential operator with coefficients supported in $M
\setminus X_{R_0} \setminus \neigh (\gamma )$.  We calculate
\ben
\label{A-op-cutoff}
\| A(x,hD)  u \| \leq \| A(x,hD) \phi^w  u \| + \| A (x, hD) (1 -
\phi)^w  u \|.
\een
Now
\be
\| A(x,hD) \phi^w  u \| \leq |E + \alpha| \| \nabla \chi \|_{L^\infty} \| \phi_2^w  u \|
\ee
for $\phi_2 \in \Ci_c (T^*M)$, $0 \leq \phi_2 \leq 1$, a microlocal cutoff supported away from
$\gamma$.  From Lemma \ref{wf-lemma-2}, we have
\be
\lefteqn{ |E + \alpha| \| \nabla \chi \|_{L^\infty} \| \phi_2^w  u \| } \\
& \leq & C R_7 h^{-1} \| Q(z)
u \| +   C \| \nabla \chi \|_{L^\infty} \| \psi_1
 u \| + \O(h) \| \tchi u \|.
\ee
For
the second term in \eqref{A-op-cutoff} choose $\chi_2 \in \Ci_c(M)$
satisfying $\chi_2 \equiv 1$ on $\supp \nabla \chi$ with support in a
slightly larger set, and $| \nabla \chi_2| \leq 2 |\nabla \chi|$.  We calculate:
\ben
\| A (x, hD) (1 -
\phi)^w u  \|^2 & = &  \lll A(x,hD)^* A(x, hD) \chi_2 (1 - \phi)^w u  ,
\chi_2 (1 -
\phi)^w u  \rrr \nonumber \\
& \leq & \half \| A^* A \chi_2 (1 - \phi)^w u  \|^2 + \half \| \chi_2 (1 - \phi)^w u 
\|^2 \nonumber \\
& \leq & \frac{C}{ \alpha} \|  (P(h) - z) \chi_2 (1 - \phi)^w u  \|^2
\nonumber 
\\
& \leq & \frac{C}{ \alpha} R_7 \| Q(z) u  \|^2 + C h^2 \|\nabla \chi \|_{L^\infty}  \| \psi_1 u  \|^2 +
\O(h^4) \|\tchi u  \|^2, \nonumber
\een
where we have again used Lemma \ref{wf-lemma-2}, \eqref{P-phi-comm} and
the fact that $P(h) -z$ is a second order elliptic semiclassical
differential operator on $\supp (1 - \phi)^w$.

We have shown
\ben
\lefteqn{ h^{-1} \left\| [ P(h), \chi ] u \right\| } \nonumber \\ && \leq CR_7 h^{-1} \| Q(z)   u \| + 
C \|\nabla \chi \|_{L^\infty} \| \psi_1  u \| + \O(h) \| \tchi u \|. 
\label{Q-101}
\een
We now use the special structure of $Q(z)$ to absorb the error terms.
To do this, choose $C_1>0$ sufficiently large that 
\be
(C_1 a - c_0 / \log(1/h)) \psi_1 \geq (C_1 a - c_0 ) \psi_1 \geq  c_0
\psi_1 / 2,
\ee
and recall $\psi_1 = \psi_2^2$.
Then
\be
\|\psi_1 u \| \leq \| \psi_2 u \|, 
\ee
and for any $\eta>0$,
\be
\half c_0 h \int  \left| \psi_2 u \right|^2 dx & \leq & h \int \left(
  C_1 a + \Im z /h \right) u \overline{ \psi_1 u } dx \\ 
& \leq & - \Im \int Q(z) u \overline{ \psi_1 u } dx \\ 
& \leq & (4 \eta h )^{-1} \| Q(z) u \|^2 + \eta h \| \psi_1 u \|^2.
\ee

Combining the last two inequalities yields
\be
\| \psi_1 u \|^2  \leq  C (4 \eta h^2 )^{-1} \| Q(z) u \|^2 + C \eta  \| \psi_1 u \|^2,
\ee
which, for sufficiently small $\eta>0$ independent of $h$, gives
\be
\| \psi_1 u \|^2 \leq C (4 \eta h^2 )^{-1} \| Q(z) u \|^2.
\ee
Plugging into \eqref{Q-101} gives
\be
h^{-1} \left\| [P(h), \chi ] u \right\| \leq C h^{-1} R_7 \| Q(z) u \| +
\O(h) \| \tchi u \|.
\ee
\end{proof}
\vspace{.25cm}
\noindent {\bf Step 4: $(P-z)^{-1}$ and $\rho_s$.}

We have shown 
\ben
\label{Q-lower-bd}
\| \psi u \|   \leq \frac{C \log (1/h)}{h} R_7 \|Q(z) u \| + \O(h)
\|\tchi u \|,
\een
but we have yet to replace $Q(z)$ in the estimate with $P(h)-z$ and add
the weights $\rho_s$.  Recall we have assumed $\rho_s \equiv 1$ on
$\supp \tchi$, and we have yet to determine the $R_j$s.  Then
\be
\| \rho_s u \| \leq \| \psi u \| + \| \rho_s (1-\psi) u \|.
\ee
Recall there is $\tP(h)$ which agrees with $P(h)$ on $\supp (1 -
\psi)$, and the principal symbol, $\tp$, of $\tP(h)$ has globally
non-trapping classical flow.  Applying \cite[Theorem 1.1]{CPV}, we get
\be
\| \rho_s (1 - \psi) u \| & \leq & C h^{-1} \| \rho_{-s}(\tP(h) -z)(1 -
\psi) u \| \\
& = & C h^{-1} \| \rho_{-s}(P(h) -z)(1 -
\psi) u \|.
\ee


Thus
\ben
\label{P-est-10}
\lefteqn {C \|\rho_{-s} (P(h)-z) u \|^2 } \\
 & \geq & C_2 \left( \| \rho_{-s} (1 - \psi) (P(h)-z) u \|^2 + \|
\psi (P(h)-z) u \|^2 \right) \\
& \geq & \| \rho_{-s} (P(h)-z) (1-\psi) u \|^2 + \| (P(h)-z) \psi u \|^2
\nonumber \\
&& + 2\|[P(h),\psi]  u \|^2 
\nonumber \\
&& - 2  \|[P(h),\psi]  u \|\left( \|  (P(h)-z) \psi u \| + \|
  \rho_{-s} (P(h)-z) (1-\psi) u\| \right). \nonumber 
\een
Applying the Cauchy-Schwarz inequality to the last term on the right
hand side yields
\be
\lefteqn{ 2 \|[P(h),\psi]  u \| \left( \|  (P(h)-z) \psi u \| + \|
  \rho_{-s} (P(h)-z) (1-\psi) u\| \right) } \\
 & \leq & 2 \|[P(h),\psi]  u \| ( \|  \psi (P(h)-z)  u \| + \|
\rho_{-s} (1 - \psi) (P(h)-z) u \| + 2\| [P(h),\psi]  u \|) \\ 
& \leq & 2 \left( 3 \| [P(h),\psi]  u \|^2 + \half \|  \psi (P(h)-z)
  u \|^2 + \half \|
\rho_{-s} (1 - \psi) (P(h)-z) u \|^2 \right).
\ee
Plugging into \eqref{P-est-10}, we have
\ben
C_3 \|\rho_{-s}(P(h)-z) u \|^2 & \geq & \| \rho_{-s} (P(h)-z) (1 - \psi) u \|^2 + \| (P(h)-z)
\psi u \|^2 \nonumber \\
&& - 4\|[P(h),\psi]  u \|^2  . \label{P-est-11}
\een
Applying Lemma \ref{wf-lemma-2} to the last term in \eqref{P-est-11} with $A
= \rho_s (1 - \psi)$, we get 
\be
\lefteqn{ \|[P(h),\psi]  u \| } \\
&& \leq C' R_7 \|(P(h)-z) u \| + \frac{C}{(R_4-R_3)^\half} h \|
\rho_s (1 - \psi) u \| + \O(h^2) \| \tchi u \|,
\ee
with $C$ here independent of $R_3$, $R_4$, and $h$.  Hence
\be
C_3 \|\rho_{-s} (P(h)-z) u \|^2 & \geq & \| \rho_{-s} (P(h)-z) (1 - \psi) u \|^2 + \| (P(h)-z)
\psi u \|^2  \\
&& - 4 \frac{C^2}{R_4 - R_3}  h^2 \| \rho_s (1 - \psi) u \|^2 - \O(h^4) \| \tchi u \|^2
\\
& \geq & \frac{h^2}{C^2} \| \rho_s (1 - \psi) u \|^2 + \frac{h^2}{C^2 \log^2
  (1/h)} \| \psi u \| \\
&& - 4 \frac{C^2}{R_4 - R_3}  h^2 \| \rho_s (1 - \psi) u \|^2 - \O(h^4) \| \tchi u \|^2 \\
& \geq & \frac{h^2}{C_4 \log^2
  (1/h)} \| \rho_s u \|^2,
\ee
as long as $R_4 - R_3>0$ is sufficiently large but fixed.  Fixing the
other $R_j$s appropriately gives \eqref{res-est}.

\end{proof}

Theorem \ref{res-thm-tau} now follows immediately from Corollary  
\ref{res-cor}.
\qed


\section{Proof of Theorems \ref{Sch-thm} and \ref{Wave-thm}}
\subsection{Proof of Theorem \ref{Sch-thm}}
In this section we show how to use the estimate \eqref{res-est} to
prove Theorem \ref{Sch-thm}.  This is an adaptation of the similar proof in \cite{Bur2}, in
the case $M$ is Euclidean space with several convex bodies removed and
compactly supported weights.

Let $\rho_s$ satisfy \eqref{rho-def}, let $\mu = \tau \pm i \epsilon$, and suppose $u$ and $f$ satisfy
\ben
\label{u-f-eqn}
 (\Delta_g  -V + \mu) u = \rho_s f.
\een
We multiply by $\rho_s^2 \bar u$ and integrate:
\be
&& \int \rho_s^2 \bar u \Delta u  + \int
(\mu -V) \rho_s^2 |u|^2 =  \int \rho_s^3 f \bar u \\
&& \implies - \int \rho_s^2 | \nabla u |^2 +  \int \mu  \rho_s^2
|u|^2 - \int ( \nabla u, \nabla (\rho_s^2)) \bar u = \int \rho_s^3 f
\bar u 
\ee
which implies
\be
\int \rho_s^2 | \nabla u |^2  & \leq & ( |\tau| + C) \int
\rho_s^2 |u|^2 + \beta \int | \nabla u |^2 | \nabla ( \rho_s^2)|^2
\rho_{-s}^2 \\
&& \quad \quad + ( 4 \beta)^{-1}  \int \rho_s^2 |u|^2 + | \int \rho_s^3
f \bar u |
\ee
for any $\beta>0$, since $V$ is bounded.  We observe 
\be
|\nabla ( \rho_s^2)|  \leq  C\lll x \rrr^{-2s-1}
\ee
for large $|x|$, and hence
\be
|\nabla ( \rho_s^2) |^2 \rho_{-s}^2 \leq C \rho_s^2
\ee
for large $|x|$.  This combined with $\rho_s^2 \leq C \rho_s$ implies
\ben
\label{u-der-est}
\int \rho_s^2 | \nabla u |^2 \leq ( | \tau | + C) \| \rho_s u
\|_{L^2}^2 + \| \rho_s  f \|_{L^2}^2.
\een
Now \eqref{u-f-eqn} implies
\be
( |\tau| + C)^{1/2} \| \rho_s u \|_{L^2} & \leq & ( |\tau| + C )^{1/2} \| \rho_s
( \Delta_g -V + \mu)^{-1} \rho_s f \|_{L^2} \\
& \leq & C \log(2 +
|\tau|) \| f \|_{L^2},
\ee
which combined with \eqref{u-f-eqn} gives
\be
\| \rho_s u \|_{H^1}^2 \leq C\int \rho_s^2 | \nabla u |^2 + C \int \| \rho_s u \|_{L^2}^2 \leq C \log(2 +
|\tau|) \| f \|_{L^2}^2.
\ee

This combined with the standard interpolation arguments gives the
following lemma.
\begin{lemma}
\label{H1-lemma}
With the notation and assumptions above, we have
\be
\| \rho_s (-\Delta_g  +V - (\tau \pm i \epsilon))^{-1} \rho_s \|_{L^2 \to H^1} \leq C_\epsilon \log (2 +
|\tau|)
\ee
and for every $\delta >0$, $r \in [-1,1]$, 
\be
\| \rho_s (-\Delta_g  +V - (\tau \pm i \epsilon))^{-1} \rho_s \|_{H^r
  \to H^{1 + r - \delta}} \leq C_{\epsilon, \delta}.
\ee
\end{lemma}

Now let $A$ be the operator
\be
A u_0 = \rho_s e^{-itP} u_0,
\ee
acting on $L^2(M)$.  We want to show
\be
A : L^2(M) \to L^2([0,T]; H^{\half - \epsilon}(M))
\ee 
is bounded.  We use the standard argument from \cite{BGT}.  That is, by duality, this is equivalent to the adjoint $A^*$ being
bounded 
\be
A^* : L^2([0,T]; H^{-\half + \epsilon}(M)) \to L^2(M),
\ee
which is equivalent to the boundedness of 
\be
A A^* : L^2([0,T]; H^{-\half + \epsilon}(M)) \to L^2([0,T]; H^{\half -
  \epsilon}(M)).
\ee
The definition of $A$ gives
\be
A^* f = \int_0^T  e^{i \tau P} \rho_s f( \tau ) d \tau 
\ee
so
\be
A A^* f(t) = \int_0^T \rho_s e^{-i  (t - \tau) P} \rho_s f( \tau) d
\tau.
\ee

We show $A A^*$ is bounded.  Let $u$ be defined by 
\be
u(x,t) = \int_0^T e^{-i  (t - \tau) P} \rho_s f( \tau) d
\tau.
\ee
Since we are only interested in the
time interval $[0,T]$, we extend $f$ to be $0$ for $t \notin [0,T]$.
We write
\be
A A^* f(t) & = & \int_0^t \rho_s e^{-i  (t - \tau) P} \rho_s f( \tau) d
\tau + \int_t^T \rho_s e^{-i  (t - \tau) P} \rho_s f( \tau) d
\tau \\
& = :& \rho_s u_1(t) + \rho_s u_2(t),
\ee
and calculate
\ben
\label{uj-eqn}
(D_t + P) u_j = (-1)^j i \rho_s f.
\een
Thus boundedness of $A A^*$ will follow if we prove $u$ satisfying
\eqref{uj-eqn} satisfies
\be
\| \rho_s u \|_{L^2([0,T];H^{\half - \epsilon})} \leq \| f
\|_{L^2([0,T];H^{-\half + \epsilon})}.
\ee
Replacing $\pm i f$ with $f$ in equation \eqref{uj-eqn} and taking the Fourier transform in time results in the following
equation for $\hat u$ and $\hat f$:
\ben
\label{ft-u-eqn}
(-z +P ) \hat u(z, \cdot ) = \rho_s \hat f (z, \cdot).
\een
Since $f(t, \cdot)$ is supported only in $[0,T]$,
$\hat f(z, \cdot)$ and $\hat u(z, \cdot)$ 
are holomorphic, bounded, and satisfy \eqref{ft-u-eqn} in $\{ \Im z <0 \}$.  Let $z =
\tau - i \eta$, $\eta >0$ sufficiently small.  We apply Lemma
\ref{H1-lemma} to get
\be
\| \rho_s \hat u(z, \cdot) \|_{H^{\half - \epsilon}(M)} \leq C \|
 \hat f(z, \cdot) \|_{H^{-\half + \epsilon}(M)},
\ee
for $\epsilon>0$.  Thus
\be
\| \rho_s  u \|_{L^2( [0,T] ; H^{\half - \epsilon}(M))} & \leq &
e^{\eta T} \| e^{-\eta t} \rho_s u(t) \|_{L^2( [0,T] ;
  H^{\half - \epsilon}(M))} \\
& \leq & Ce^{\eta T} \| \rho_s \hat{u}(\tau - i \eta) \|_{L^2(
  \reals ; H^{\half - \epsilon}(M))} \\
& \leq & Ce^{\eta T} \| \hat f( \tau - i \eta) \|_{L^2( \reals ; H^{-\half +
    \epsilon}(M))} \\
& \leq & C e^{\eta T} \|e^{- \eta t} f(t) \|_{L^2([0,T];  H^{-\half +
    \epsilon}(M))} \\
& \leq & C e^{\eta T} \| f(t) \|_{L^2([0,T];  H^{-\half +
    \epsilon}(M))}.
\ee
Hence 
\be
\int_0^T \|\rho_s u \|_{H^{\half - \epsilon}(M)}^2 dt \leq C e^{ \eta
  T} \int_0^T \| f \|_{H^{-\half + \epsilon}(M)}^2 dt,
\ee
or $A A^*$ is bounded.  Thus $A$ is bounded and Theorem \ref{Sch-thm}
is proved.
\qed

\begin{remark}
\label{rem-2}
If the estimate \eqref{res-tau-est} is uniform in the lower
half-plane, then the preceding calculation can be made including
taking the limit $\eta \to 0$, in which case we get the global in time
local smoothing estimate \eqref{global-sm}
\end{remark}

The following Lemma uses interpolation to replace the $H^{1/2- \epsilon}$ norm
on the left hand side of \eqref{Sch-thm-est} with $H^{1/2}$, and will be
of use in \S \ref{strichartz}.

\begin{lemma}
\label{inter-lem-10}
Suppose $(M,g)$ and $V$ satisfy the assumptions of Theorem
\ref{Sch-thm}.  For each $\delta>0$ there is a constant $C>0$ such
that
\ben
\label{eqn-1002}
\int_0^T \left\| \rho_s e^{it (\Delta_g - V(x)) } u_0 \right\|_{H^{1/2}(M)}^2 dt  \leq C \| u_0 \|_{H^\delta(M)}^2.
\een
\end{lemma}

\begin{proof}
We first calculate
\be
\| \rho_s e^{itP} u_0 \|_{L^2_T H^{1}}^2 & \leq & C \int_0^T \int_M
\left| \rho_s (P+1)e^{itP} u_0 \rho_s e^{-itP} \overline{u_0} \right|
dx dt \\
&& + 2 \int_0^T \int_M | \nabla \rho_s| \left| \nabla e^{itP}u_0 \right|
\left| \rho_s e^{itP} u_0 \right| dx dt \\
&& + \int_0^T \int_M |P( \rho_s) | \left| e^{itP} u_0 \right| \left| \rho_s
  e^{itP} u_0 \right| dx dt.
\ee
Using 
\be
| P( \rho_s)| \leq C | \nabla \rho_s | \leq C' |\rho_s|
\ee
and applying the Cauchy-Schwarz inequality yields
\be
\|\rho_s e^{itP} u_0 \|_{L^2_T H^1}^2 \leq C \| u_0 \|_{H^2}.
\ee
Thus we have a linear operator bounded between complex interpolation
spaces:
\be
\rho_s e^{itP} & : & L^2  \to L^2_T H^{1/2 - \epsilon}, \\
& : & H^2 \to L^2_T H^1.
\ee
Choosing $\epsilon = \delta/4$ we have
\be
\| \rho_s e^{itP} u_0 \|_{L^2_T H^{1/2}} & \leq & C \| u_0 \|_{H^{2
    \epsilon/ (1/2 + \epsilon)}} \\
& \leq & C \| u_0 \|_{H^\delta}.
\ee
\end{proof}

\subsection{Proof of Theorem \ref{Wave-thm}}
For the proof of Theorem \ref{Wave-thm}, we apply \cite[Theorem
3]{Ch3}, which is a generalization of \cite[Th\'eor\`eme 3]{Bur1a}.
That is, we set
\be
B = \left( \begin{array}{cc} 0 & -i \id \\ -i \Delta & 0 \end{array}
\right),
\ee
acting on the Hilbert space $H = H^1(X) \times L^2(X)$.  The commutator
$[\psi, B]$ is bounded on $H$, so if $\psi_2 \in \Ci(X)$ satisfies
\eqref{psi-def} and $|[\psi, -\Delta_g]| \leq \psi_2$, we have
\be
\| \psi e^{itB} \psi \|_{\Dom (B^2) \to  H} & = & \| \psi e^{itB} \psi
(1 -iB)^{-2}\|_{H \to H} \\
& \leq & C \| \psi e^{itB} 
(1 -iB)^{-2} \psi_2 \|_{H \to H}.
\ee
From \cite[Theorem 3]{Ch3} we then gather
\be
\| \psi e^{itB} \psi \|_{\Dom (B^2) \to  H} \leq C e^{-t^{1/2} /C}.
\ee
The spaces $H^{1+s} \times H^s$ are complex interpolation spaces, so
together with the trivial estimate
\be
\| \psi e^{itB} \psi \|_{H \to  H} \leq C,
\ee
we conclude that for any $\epsilon>0$,
\be
E_\psi (t) \leq C_\epsilon e^{-\epsilon t^{1/2} /C} \left( \|u_0
  \|_{H^{1+ \epsilon} }^2 + \|u_1 \|_{H^\epsilon}^2 \right).
\ee

\qed

\section{Strichartz-type Inequalities}
\label{strichartz}
In this section we prove several families of Strichartz-type
inequalities and prove Proposition \ref{nls-lwp}.  The statements and
proofs are mostly adaptations of
similar inequalities 
in \cite{BGT}, so we leave out the proofs of these in the interest of space.

As in the statement of Proposition \ref{nls-lwp}, we assume $M$ is
asymptotically conic as defined in \cite{HTW} and $V \in \Ci_c(M)$, $V
\geq 0$.  The manifold $M$ admits the Sobolev embeddings recorded in the following proposition.  For
our notation, let
\be
W^{m,p}(M)\,\,\, (\text{resp. } W_0^{m,p}(M)), \,\, m \in \mathbb{N}
\ee
be the completion of $\Ci(M)$ (resp. $\Ci_c(M)$) with respect to the norm
\be
\|f\|_{W^{m,p}}^p = \sum_{|\alpha| \leq m} \|D^\alpha f \|_{L^p}^p.
\ee
We define $W^{s,p}(M)$ and $W_0^{s,p}(M)$ for non-integer $s$ by
interpolation.  We use
the convention
\be
H^s(M):= W^{s,2}(M), \text{ and } H_0^s(M):= W_0^{s,2}(M).
\ee
Let $H_D^1(M)$ denote the domain of $(1- \Delta_g)^\half$ with
Dirichlet boundary conditions if $\partial M \neq \emptyset$, so that $H_D^1(M)
= H_0^1(M)$, and write $H_D^s(M)$ for the domain of $(1 -
\Delta_g)^{s/2}$ with Dirichlet boundary conditions.  Since $V \in \Ci_c(M)$, we may replace $(1-
\Delta_g)^\half$ with $(1 + P)^\half$ in the definitions, where $P = - \Delta_g + V(x)$.  This results in equivalent
Sobolev spaces with the addition $[P, (1+P)^\half]=0$.

\begin{proposition}
\label{sobolev-prop}
We have the following continuous Sobolev embeddings:
\be
&& \text{(i) } H_D^1(M) \subset L^p(M), \,\,\, 2 \leq p \leq
\frac{2n}{n-2}, \,\, \text{or } p < \infty \text{ for } n=2, \\
&& \text{(ii) } H_D^s(M) \subset L^p(M), \,\,\, \frac{1}{2} =
\frac{s}{n} + \frac{1}{p}, \,\, s \in [0,1), \\
&& \text{(iii) } H_D^{s+1}(M) \subset W^{1,p}(M), \,\,\, \frac{1}{2} =
\frac{s}{n} + \frac{1}{p}, \,\, s \in [0,1), \\
&& \text{(iv) } W_0^{1,p}(M) \subset L^q(M), \,\,\, \frac{1}{p} =
\frac{1}{n} + \frac{1}{q}, \,\, 1 \leq p < q < + \infty, \\
&& \text{(v) } W_0^{s,p}(M) \subset L^\infty (M), \,\,\, s >
\frac{n}{p}, \,\, p \geq 1 \\
&& \text{(vi) } H_D^{s + 1/p}(M) \subset W^{s ,q}(M), \,\,\,
\frac{1}{p} + \frac{n}{q} = \frac{n}{2}, \,\, p \geq 2, s \in [0,1].
\ee
\end{proposition}

If we again let $-\Delta_0$ be the Laplace-Beltrami operator
associated to a non-trapping metric which agrees with $g$ on $X_R$,
we may apply the results of \cite{HTW} to a solution of the
Schr\"odinger equation away from the trapping region, resulting in
perfect Strichartz estimates, but we lose something from the commutator.  That is, if $\chi \in \Ci_c(M)$ is $1$
on $M \setminus (X_R \cup \supp V)$, then $w = (1 - \chi) e^{-itP}u_0$ satisfies
\ben
(D_t-\Delta_0 ) w & = & (D_t + P) w \nonumber \\
& = & [\Delta_0, \chi] e^{-itP}u_0. \label{sch-eq-comm}
\een

From Lemma \ref{inter-lem-10}, we have for any $\epsilon>0$,
\be
\|[\Delta_0, \chi] e^{-itP}u_0 \|_{L^2([0,T]) H^{-1/2}(M)} \leq
  C_{\epsilon, T} \|u_0 \|_{H^\epsilon(M)}.
\ee
The following proposition then follows from the proof of
\cite[Proposition 2.10]{BGT}.

\begin{proposition}
\label{str-prop-1}
For every $0 < T \leq 1$, $\delta>0$, and each $\chi \in \Ci_c(M)$ satisfying $\chi
\equiv 1$ near $M \setminus (X_R \cup \supp V)$ , there is a constant $C>0$ such that
\ben
\label{str-1}
\| (1-\chi) u \|_{L^p([0,T]) W^{s- \delta,q}(M)} \leq C \| u_0 \|_{H_D^s(M)},
\een
where $u = e^{-it P}u_0$, $s \in [0,1]$, and $(p,q)$, $p >2$
satisfy
\be
\frac{2}{p} + \frac{n}{q} = \frac{n}{2}.
\ee
\end{proposition}

\begin{remark}
In the sequel, wherever unambiguous, we will write
\be
L_T^pW^{s,q}:= L^p([0,T]) W^{s,q}(M)
\ee
and
\be
H_D^s: = H_D^s(M).
\ee
\end{remark}



\begin{proposition}
\label{str-loss-1}
Let $u(t) = e^{-it P} u_0$.  For every $0 < T \leq 1$ and $\epsilon>0$, there is a constant $C>0$
such that
\ben
\label{str-2}
\|u \|_{L_T^p W^{s,q}} \leq C \| u_0 \|_{H_D^{s + 1/p + \epsilon}},
\een
where $s \in [0,1]$ and $(p,q)$, $p>2$ satisfy
\be
\frac{2}{p} + \frac{n}{q} = \frac{n}{2}.
\ee
\end{proposition}
\begin{remark}
Proposition \ref{str-loss-1} represents the 
Strichartz estimates obtained by Burq-G\'erard-Tzvetkov \cite{BGT} in
the case of non-trapping exterior domains with an $\epsilon>0$ loss due to the presence of
the trapped orbit $\gamma$.  Observe that \eqref{str-loss-1} is weaker
than the standard Euclidean Strichartz estimates in two ways, the loss
of $1/p$ derivatives from using Sobolev embeddings and the loss of
$\epsilon$ derivatives from $\gamma$.  When $\partial M = \emptyset$,
we get the improvement given in Proposition \ref{str-prop-5}.
\end{remark}

\begin{proposition}
\label{str-prop-3}
Let $u = e^{-it P} u_0$ and 
\be
v = \int_0^t e^{-i(t - \tau) P} f( \tau ) d \tau.
\ee
Then for each $0 <T \leq 1$ and each $\delta>0$, there exists $C >0$
such that
\ben
\label{str-3}
\| u \|_{L^p_T W^{s- \delta, q}} \leq C \| u_0 \|_{H_D^s}
\een
and
\ben
\label{str-4}
\| v \|_{L_T^p W^{s - \delta, q}} \leq C \| f \|_{L_T^1 H_D^s},
\een
where $s\in [0,1]$ and $(p,q)$, $p>2$ satisfy
\ben
\label{pq-eqn-2}
\frac{1}{p} + \frac{n}{q} = \frac{n}{2}.
\een
\end{proposition}

\begin{remark}
Proposition \ref{str-prop-3} is much weaker than the estimate
suggested by scaling in Euclidean space, and as remarked in
\cite{BGT}, is probably not optimal.  We expect the $\delta>0$ loss to
always hold due to the presence of $\gamma$, but the Euclidean scaling
suggests the optimal estimate would replace $1/p$ in \eqref{pq-eqn-2}
with $2/p$ (see Proposition \ref{str-prop-5}).
\end{remark}

\begin{proposition}
\label{str-prop-4}
Let 
\be
v(t) = \int_0^t e^{-i(t - \tau) P} f( \tau) d \tau.
\ee
For each $0 < T \leq 1$ and each $\delta>0$, there is a constant $C>0$
such that
\ben
\label{str-5}
\| v \|_{L_T^p W^{s - \delta, q}} \leq C \| f \|_{L_T^{p'} W^{s, q'}},
\een
where $p', q'$, $p' \in [1,2)$ are the duals of $p$ and $q$ satisfying \eqref{pq-eqn-2},
respectively, and satisfy
\be
\frac{1}{p'} + \frac{n}{q'} = \frac{n}{2} + 1.
\ee
\end{proposition}

The next proposition is an improvement of Proposition \ref{str-prop-3} in the
case $\partial M = \emptyset$.  
\begin{proposition}
\label{str-prop-5}
Suppose $(M,g)$ and $V$ satisfy the assumptions of Proposition
\ref{nls-lwp}, $u = e^{-it P} u_0$, 
\be
v = \int_0^t e^{-i(t - \tau) P} f( \tau ) d \tau,
\ee
and in addition $\partial M = \emptyset$.  Then for each $0 <T \leq 1$ and each $\delta>0$, we have the estimates
\eqref{str-3} and \eqref{str-4} for $s\in [0,1]$, where now $(p,q)$,
$p>2$ satisfy the Euclidean scaling 
\ben
\label{pq-eqn-3}
\frac{2}{p} + \frac{n}{q} = \frac{n}{2}.
\een
\end{proposition}
\begin{proof}
The idea of the proof is to use Proposition \ref{str-prop-1} to reduce
the statement to a local question near the trapped orbit.  Then we use
a partition of unity and 
the local WKB construction from \cite{BGT1} to get local in time
Strichartz estimates for time on the scale of inverse frequency.  We
then sum up over frequencies and apply the local smoothing estimate to
prove the Proposition.  We remark this would also follow from
\cite[Theorem 4]{StTa} and the local smoothing in Theorem \ref{Sch-thm}.

Let $\chi$ be as in Proposition \ref{str-prop-1} and choose $\psi \in
\Ci_c(\reals)$, $\psi \equiv 1$ near $1$ and satisfying
\be
1 \leq \sum_{k \geq 0} \psi(  r/k ) \leq 2 \text{ for } r \geq 0.
\ee
Choose also $\phi \in \Ci_c( \reals)$, $\phi \equiv 1$ on $[-c_0,
c_0]$, $\supp \phi \subset [-2c_0, 2c_0]$ for $c_0>0$ small.  Let 
\be
w_h = \phi(t/h) \chi(x) \psi( -h^2 \Delta_g + h^2 V(x)) u,
\ee
which satisfies the equation
\be
\left\{ \begin{array}{l}
(i \partial_t + \Delta_g - V(x))w_h = \phi [\Delta_g, \chi ] \psi u +
i\frac{1}{h} \phi' \chi \psi u \\ w_h(x,0) = \chi \psi u_0. \end{array}
\right.
\ee
Since $\phi$ localizes to a timescale of size $h$, the semiclassical
local WKB construction in \cite{BGT1} gives
\be
\| w_h \|_{L^p L^q} \leq C \left\| \phi [\Delta_g, \chi ] \psi u +
i\frac{1}{h} \phi' \chi \psi u \right\|_{L^1 L^2},
\ee
with $(p,q)$, $p>2$ satisfying \eqref{pq-eqn-3}.

Choose $\tphi \in \Ci_c( \reals)$ and $\tchi \in \Ci_c(M)$ satisfying
$\tphi \equiv 1$ on $\supp \phi$ and $\tchi \equiv 1$ on $\supp
\chi$.  Applying H\"older's inequality in time and exchanging one
derivative for $h^{-1}$ on the support of $\psi$ yields
\be
\| w_h \|_{L^p L^q} \leq C h^{1/2}  \left\|h^{-1} \tphi \tchi \psi u
\right\|_{L^2 L^2}.
\ee
Exchanging one half derivative with $h^{-1/2}$ we obtain
\be
\| w_h \|_{L^p L^q} \leq C  \left\| \tphi \tchi \psi u
\right\|_{L^2 H^{1/2}},
\ee
and summing in $h=1/k$ we get
\be
\| \chi u \|_{L^p L^q} \leq C \| \tchi u \|_{L^2 H^{1/2}},
\ee
which after a time truncation and an application of Lemma \ref{inter-lem-10}
proves the Proposition for $u$.  Finally, an application of the
Christ-Kiselev lemma \cite{CK} proves the proposition for $v$.
\end{proof}


\begin{proof}[Proof of Proposition \ref{nls-lwp}]
The proof of Proposition \ref{nls-lwp} is a slight modification of the
proof of Proposition 3.1 in \cite{BGT1}, but we include it here in the
interest of completeness.  First we assume $\partial M \neq
\emptyset$.  
Fix $s$ satisfying \ref{s-cond} and choose $p > \max \{ 2\beta-2, 2 \}$
satisfying
\be
s > \frac{n}{2} - \frac{1}{p} + \delta \geq \frac{n}{2} - \frac{1}{\max
  \{ 2\beta-2, 2 \}}
\ee
for some $\delta>0$.  Set $\sigma = s - \delta$ and 
\be
Y_T = C ( [-T,T]; H_D^s(M)) \cap L^p([-T,T]; W_0^{\sigma,q}(M))
\ee
for 
\be
\frac{1}{p} + \frac{n}{q} = \frac{n}{2},
\ee
equipped with the norm
\be
\| u \|_{Y_T} = \max_{|t| \leq T} \| u(t)\|_{H^s_D(M)} + \|u\|_{L^p_T
  W^{\sigma,q}}.
\ee
Let $\Phi$ be the nonlinear functional
\be
\Phi ( u ) = e^{-it P} u_0 - i \int_0^t e^{-i(t-\tau) P}
F(u(\tau)) d \tau.
\ee
If we can show that $\Phi: Y_T \to Y_T$ and is a contraction on a ball
in $Y_T$ centered at $0$ for sufficiently small $T>0$, this will
prove the first assertion of the Proposition, along with the Sobolev
embedding
\ben
\label{sob-emb-1a}
W^{\sigma,q}_0(M) \subset L^\infty(M),
\een
since $\sigma >n/q$.  From Proposition \ref{str-prop-3}, we bound the
$W^{\sigma}$ part of the $Y_T$ norm by the $H^s_D$ norm, giving
\be
\| \Phi(u) \|_{Y_T} & \leq & C \left( \| u_0 \|_{H^s_D} + \int_{-T}^T
  \| F(u (\tau)) \|_{H^s_D} d \tau \right)\\
& \leq & C \left( \| u_0 \|_{H^s_D} + \int_{-T}^T \|(1 +
  |u(\tau)|)\|_{L^\infty}^{2\beta-2} ) \|u( \tau ) \|_{H^s_D} d \tau
\right),
\ee
where the last inequality follows by our assumptions on the structure
of $F$.  
Applying H\"older's inequality in time with $\tilde p = p/ (2\beta-2)$
and
\be
\frac{1}{\tilde q} + \frac{1}{ \tilde p} = 1
\ee
gives
\be
\| \Phi(u) \|_{Y_T}  \leq C \left( \| u_0 \|_{H^s_D} + T^\gamma  \|
  u\|_{L^\infty_T H^s_D}\|( 1 + |u| )\|_{L^p_T L^\infty}^{2\beta-2})
\right)
\ee
where $\gamma = 1/ \tilde q >0$.  Thus
\be
\| \Phi(u) \|_{Y_T} \leq  C \left( \| u_0 \|_{H^s_D} +
  T^\gamma(\|u\|_{Y_T} + \|u\|_{Y_T}^{ 2\beta} ) \right).
\ee
Similarly, we have for $u, v \in Y_T$,
\ben
\label{phi-ineq-2} \lefteqn{ \| \Phi(u) - \Phi(v) \|_{Y_T}  \leq } \\
 & \leq &  C T^\gamma \| u - v
\|_{L^\infty_T H^s_D} \|( 1 + | u|) \|_{L^p_T L^\infty}^{2\beta-2} +
\|(1 + |v|) 
\|_{L^p_T L^\infty}^{2\beta-2}) \\
& \leq & C T^\gamma \| u - v \|_{Y_T} \|( 1 + |u|) \|_{Y_T}^{2\beta-2} +
\|(1+ |v|) \|_{Y_T}^{2\beta-2}), \nonumber
\een
which is a contraction for sufficiently small $T$.  This concludes the
proof of the first assertion in the Proposition.  

To get the second assertion, we observe from \ref{phi-ineq-2} and the
definition of $Y_T$, if $u$ and $v$ are two solutions to \eqref{nls}
with initial data $u_0$ and $u_1$ respectively, so 
\be
\widetilde{\Phi}(v) =  e^{-it P} u_1 - i \int_0^t e^{-i(t-\tau) P}
F(v(\tau)) d \tau,
\ee
 we have
\be
\lefteqn{\max_{|t| \leq T } \| u(t) - v(t) \|_{H^s_D} } \\
& = & \max_{|t| \leq T } \| \Phi(u)(t) - \widetilde{\Phi} (v)(t) \|_{H^s_D} \\
& \leq C &
\Bigg( \| u_0 - u_1 \|_{H^s_D} \\
&& \quad + T^\gamma \max_{|t| \leq T } \| u(t) -
  v(t) \|_{H^s_D} \|( 1 + | u|) \|_{L^p_T L^\infty}^{2\beta-2} + \| (1
  + |v|)
  \|_{L^p_T L^\infty}^{2\beta-2} ) \Bigg),
\ee
which, for $T>0$ sufficiently small gives the Lipschitz continuity.

In the case $\partial M = \emptyset$, we have the improved Strichartz
estimates given in Proposition \ref{str-prop-5}.  Hence we can take
$s$ and $p$ satisfying $p > \max \{ 2 \beta -2, 2 \}$ and 
\be
s > \frac{n}{2} - \frac{2}{p} + \delta \geq \frac{n}{2} - \frac{2}{\max
  \{ 2\beta-2, 2 \}}
\ee
for $\delta>0$.  Then $\sigma = s- \delta > q/n$ and the preceding
argument holds with these modifications.
\end{proof}

The proof of Corollary \ref{gwp-cor} now follows from the standard
global well-posedness arguments from, for example, \cite[Chapter 6]{Caz}.

\qed


\end{document}